\documentclass{amsart} 
\usepackage{amssymb,latexsym}

\numberwithin{equation}{section}

\newtheorem{theorem}{Theorem}[section]

\newtheorem{lemma}[theorem]{Lemma}

\def\proof{\smallskip\noindent {\bf Proof. }}
\def\endproof{\hfill$\square$\medskip}

\def\ii{\mathbf{i}}
\def\jj{\mathbf{j}}

\def\wnot{w_\mathrm{o}}

\def\CC{\mathbb{C}}
\def\l{\ell}

\newcommand{\mat}[4]{\left(\!\!\begin{array}{cc}
#1 & #2 \\
#3 & #4 \\
\end{array}\!\!\right)}

\begin{document}

\title[Totally nonnegative and oscillatory elements]
{Totally nonnegative and oscillatory\\ elements in semisimple groups}

\author{Sergey Fomin}
\address{Department of Mathematics, Massachusetts Institute of
  Technology, Cambridge, Massachusetts 02139}
\email{fomin@math.mit.edu}

\author{Andrei Zelevinsky}
\address{\noindent Department of Mathematics, Northeastern University,
  Boston, MA 02115} 
\email{andrei@neu.edu} 

\thanks{The authors were supported in part 
by NSF grants \#DMS-9625511 and \#DMS-9700927.
}

\date{November 16, 1998}

\subjclass{
Primary 22E46; 
Secondary 
14M15, 
15A48, 
20F55. 
}

\keywords{Total positivity, oscillatory element, semisimple Lie group}

\begin{abstract}
We generalize the well known characterizations of 
totally nonnegative and oscillatory matrices,
due to F.~R.~Gantmacher, M.~G.~Krein, A.~Whitney, 
C.~Loewner, 
M.~Gasca, and J.~M.~Pe\~na
to the case of an arbitrary complex semisimple Lie group. 
\end{abstract}

\maketitle

\section{Introduction}

In this note, we extend some classical theorems in
the theory of total positivity to the case of an arbitrary semisimple
complex Lie group. We begin by reviewing the results we are going to
generalize.  

\medskip

Let $G = GL_n (\CC)$ or $SL_n (\CC)$.
The following theorem is due to 
C.~Loew\-ner~\cite{loewner} and A.~Whitney \cite{whitney}
(cf.\ \cite[Lemma~9.1]{karlin}). 

\begin{theorem} 
\label{th:Whitney-Loewner classical}
For a matrix $x \in G$, the following are equivalent:
\begin{itemize}
\item[(a)] 
all minors of $x$ are nonnegative real numbers;
\item[(b)]
$x$ lies in the closure of the set of matrices with positive
minors;  
\item[(c)]
$x$ belongs to the multiplicative monoid in $G$ 
generated by elementary Jacobi matrices with nonnegative matrix entries.
\end{itemize}
{\rm(Here in~{\rm(c)}, an ``elementary Jacobi matrix'' is a matrix that differs
from the identity matrix in a single entry, located either on the main
diagonal, or immediately above or below it.) }
\end{theorem}

A matrix $x \!\in\! G$ satisfying any of the
equivalent conditions (a)--(c) above is called \emph{totally nonnegative}.
Furthermore, $x$ is \emph{totally positive} if all its minors are
positive. 
Totally positive matrices are distinguished among the totally
nonnegative ones as follows.

\begin{theorem} 
\label{th:TNN-TP classical}
For a totally nonnegative matrix $x\!\in\! G$, the following are equivalent:
\begin{itemize}
\item[(d)] 
$x$ is totally positive; 
\item[(e)]
all 
solid minors of $x$ involving either $x_{1n}$ or $x_{n1}$ are
positive;  
\item[(f)]
$x$ belongs to the intersection of opposite open Bruhat cells 
$B\wnot B\cap B_-\wnot B_-\,$. 
\end{itemize} 
{\rm (Here in (e), a ``solid minor'' is a minor formed by several
  consecutive rows and as many consecutive columns. 
In (f), we denote by $B$ (resp.\ $B_-$) the subgroup
of upper-triangular (resp.\ lower-triangular) matrices, and $\wnot$ is
the permutation matrix with $1$'s on the main antidiagonal.) 
}
\end{theorem}

Part ${\rm (d)}\!\Longleftrightarrow\!{\rm (e)}$ of
Theorem~\ref{th:TNN-TP classical} is a refinement of the classical
Fekete criterion due to M.~Gasca and 
J.~M.~Pe\~na~\cite[Theorem~4.3]{gasca-pena};
the equivalence  ${\rm (e)}\!\Longleftrightarrow\!{\rm (f)}$ is
a well-known (and easy) linear-algebraic fact
(cf., e.g., \cite[Theorem~II.4.1]{gantmacher}). 

\medskip

In their pioneering study of total positivity undertaken in 1930s,
F.~R.~Gantmacher and M.~G.~Krein introduced and studied 
the intermediate class of \emph{oscillatory} matrices 
defined as follows: a matrix $x \in G$ is called 
oscillatory if $x$ is totally nonnegative while some power of $x$ is
totally positive. 
The following characterization of this class was obtained
in~\cite{GK} (see~\S II.7; cf.\ also \cite[Theorem~9.3]{karlin}).   

\begin{theorem} 
\label{th:GK classical}
For a totally nonnegative matrix $x\!\in\! G$, the following are equivalent:
\begin{itemize}
\item[(g)] 
$x$ is oscillatory; 
\item[(h)]
$x_{i,i+1} > 0$ and $x_{i+1,i} > 0$ for $i = 1, \dots, n\!-\!1$. 
\end{itemize} 
\end{theorem}

Gantmacher and Krein \cite{GK}  further showed that the definition 
of oscillatory matrices can be refined as follows.

\begin{theorem} 
\label{th:GK bound}
A totally nonnegative matrix $x\!\in\! G$ is oscillatory if and only if
$x^{n-1}$ is totally positive.
\end{theorem}

In this paper, we extend 
Theorems~\ref{th:Whitney-Loewner classical},
\ref{th:GK classical}, and~\ref{th:GK bound} to 
an arbitrary semisimple complex Lie group~$G$, 
using the notion of \emph{generalized minors} introduced
in~\cite{FZ}. 
A~generalization of Theorem~\ref{th:TNN-TP classical} follows from 
results in~\cite{lusztig} or~\cite{FZ},
and is presented below (see Theorem~\ref{th:TNN-TP general})
for the sake of completeness.

Even in the case of $SL_n$, our version of the criterion~(h) is more
general than the one given above. (Earlier in~\cite{FZ}, we gave a
family of total positivity criteria generalizing~(e).)
It should also be noted that our proofs are quite different from the 
ones in~\cite{gasca-pena, GK, karlin, loewner, whitney}.
Our main technical tools involve combinatorics of reduced words
in Weyl groups, the subdivison of a semisimple group into 
double Bruhat cells, and the ``generalized determinantal calculus"
developed 
in \cite{FZ}; in particular, the fundamental role 
is played by a generalized determinantal
identity~\cite[Theorem~1.17]{FZ}. 

The study of total positivity in reductive groups other than $GL_n$
and $SL_n$ was initiated by G.~Lusztig~\cite{lusztig}, who suggested
to use the natural generalization of~(c) as the definition of a total
nonnegative element. Our extension of the equivalence 
${\rm(a)}\Longleftrightarrow{\rm (c)}$ can be rephrased as saying that
Lusztig's definition is equivalent to the one in terms of the generalized
minors of~\cite{FZ}. 

\section{Terminology and notation}

We will use the setup of~\cite{FZ}, which is briefly reviewed below in
this section.
Proofs and further details can be found in~\cite{FZ} (see 
Sections~1.1-1.4). 

Let $G$ be a simply connected semisimple complex Lie group of rank~$r$
with a fixed pair of opposite Borel subgroups $B_-$ and~$B$;
thus $H=B_-\cap B$ is a maximal torus in~$G$. 
Let $N_-$ and $N$ be the unipotent radicals of $B_-$ and~$B$,
respectively.

Let $\alpha_1, \ldots, \alpha_r$ be the system of simple roots for
which the corresponding root subgroups are contained in~$N$. 
For every $i\in [1,r]$, let $\varphi_i: SL_2 \to G$ be 
the canonical embedding corresponding to the simple
root~$\alpha_i\,$. 

For any nonzero $t\in\CC$ and any $i\in [1,r]$, we define 
\begin{equation*}
\label{eq:x,y}
x_{\overline i} (t) = \varphi_i \mat{1}{0}{t}{1}\, ,\qquad
t^{h_i} = \varphi_i \mat{t}{0}{0}{t^{-1}}\, ,\qquad
x_i (t) = \varphi_i \mat{1}{t}{0}{1} \, .
\end{equation*}
Thus $t^{h_i}\in H$,
and $t \mapsto x_i (t)$ (resp.\ $t \mapsto x_{\overline i} (t)$) 
is a one-parameter subgroup in $N$ (resp.\ in $N_-$).

The \emph{weight lattice} $P$ can be defined as the group of
multiplicative characters of~$H$,
here written in the exponential notation:
a weight $\gamma\in P$ acts by $a \mapsto a^\gamma$. 
The lattice $P$ has a $\mathbb{Z}$-basis formed by the 
\emph{fundamental weights}  
$\omega_1, \ldots, \omega_r$ defined by 
$(t^{h_j})^{\omega_i} = t^{\delta_{ij}}$.  
 
The \emph{Weyl group} $W$ of $G$ is defined by $W = {\rm Norm}_G (H)/H$. 
The action of~$W$ on~$H$ by conjugation gives rise to the action 
of $W$ on the weight lattice $P$ given by
$a^{w (\gamma)} = (w^{-1} a w)^\gamma$
for $w \in W$, $a \in H$, $\gamma \in P$. 
The group $W$ is a Coxeter group with simple reflections
$s_1,\dots,s_r$ which can be defined by specifying their representatives
in ${\rm Norm}_G (H)$: we set 
\[
\overline {s_i} = \varphi_i \mat{0}{-1}{1}{0} \in {\rm Norm}_G (H) \, .
\]
The family  
$\{\overline {s_i}\}$ 
satisfies the braid relations in~$W$;
thus the representative $\overline w$ 
can be unambiguously defined for any 
$w \in W$ by requiring that 
$\overline {uv} = \overline {u} \cdot \overline {v}$
whenever $\l (uv) = \l (u) + \l (v)$;
here $\l(w)$ denotes the length of~$w\in W$.

A \emph{reduced word} for $w \in W$ is a sequence of indices 
$(i_1, \ldots, i_m)$ that satisfies $w = s_{i_1} \cdots s_{i_m}$
and has the shortest possible length~$m=\l(w)$. 
The set of reduced words for~$w$ will be denoted by~$R(w)$. 

As customary, $\wnot$ denotes the unique element of maximal length
in~$W$. 

We denote by $G_0=N_-HN$ the open subset of elements $x\in G$ that 
have Gaussian decomposition; this decomposition will be written as
$x = [x]_- [x]_0 [x]_+ \,$.

For $u,v \in W$ and $i \in [1,r]$, the \emph{generalized minor} 
$\Delta_{u \omega_i, v \omega_i}$ 
is the regular function on $G$ whose restriction to the open set
${\overline {u}} G_0 {\overline {v}}^{-1}$ is given by
\begin{equation*}
\Delta_{u \omega_i, v \omega_i} (x) = 
\left[{\overline {u}}^{-1}  
   x \overline v\right]_0^{\omega_i} \ . 
\end{equation*}
It can be shown that $\Delta_{u \omega_i, v \omega_i}$ depends on  
the weights $u \omega_i$ and $v \omega_i$ alone, not on the particular
choice of $u$ and~$v$. 
In the special case $G=SL_n\,$, the generalized minors are nothing but
the ordinary minors of a matrix. 

\section{Main results}


We generalize the Loewner-Whitney Theorem
(Theorem~\ref{th:Whitney-Loewner classical}) as follows. 

\begin{theorem}
\label{th:loewner general}
For an element $x \in G$, the following are equivalent:
\begin{itemize}
\item[(a)]
all generalized minors $\Delta_{\gamma, \delta}$ take nonnegative 
real values at~$x$; 
\item[(b)]
$x$ lies in the closure of the set of elements with positive
generalized minors;
\item[(c)]
$x$ lies in the multiplicative monoid generated by the elements
of the form $t^{h_i}$, $x_i (t)$, and $x_{\overline i} (t)$, with
positive~$t$.  
\end{itemize}
\end{theorem}

An element $x\in G$ satisfying any of the
equivalent conditions (a)--(c) of Theorem~\ref{th:loewner general} 
is called \emph{totally nonnegative.}
The set of all such elements is denoted by~$G_{\geq 0}\,$. 

The following generalization of Theorem~\ref{th:TNN-TP classical} is
immediate from the results in Lusztig~\cite{lusztig};
a proof based on the results in~\cite{FZ} will be given in Section~\ref{sec:TNN-TP
general} below. 

\begin{theorem} 
\label{th:TNN-TP general}
For an element $x\!\in\! G_{\geq
    0}\,$, the following are equivalent:
\begin{itemize}
\item[(d)] 
all generalized minors of $x$ are positive; 
\item[(e)]
$\Delta_{\omega_i,\wnot\omega_i}(x)>0$ and
$\Delta_{\wnot\omega_i,\omega_i}(x)>0$ for any $i\in [1,r]$; 
\item[(f)]
$x$ belongs to the intersection of open Bruhat cells 
$B\wnot B\cap B_-\wnot B_-\,$. 
\end{itemize} 
\end{theorem}

An element $x\!\in\! G$ satisfying any of the
equivalent conditions (d)--(f) of Theorem~\ref{th:TNN-TP general}
is called \emph{totally positive.}
The set of all such elements  will be denoted by~$G_{>0}\,$. 

Let us call an element $x \in G_{\geq 0}$ \emph{oscillatory} 
if for some positive integer~$m$,  
the element $x^m$ is totally positive. 
We will give equivalent reformulations of this property 
which in particular generalize the criterion (h) in 
Theorem~\ref{th:GK classical}. 
In fact, our version of criterion~(h) will be more general even in the 
special case~$G=SL_n\,$. 

Let $i$ and $j$ be two indices lying in the same connected component
of the Dynkin graph of $G$ (the case $j = i$ is not excluded). 
Let 
\[
i = i(1), i(2), \dots, i(l) = j
\] 
be the unique path from $i$ to $j$ in the Dynkin graph.
Thus $\{i(k), i(k+1)\}$ is an edge 
for $k = 1, \dots, l-1$, and all indices $i(k)$ are distinct. 
Let us denote $c(j \to i) = s_{i(2)} s_{i(3)} \cdots s_{j}$
(in particular, $c(i \to i) = e$), and set 
\begin{eqnarray*}
\begin{array}{l}
\Delta_{j \to i} = \Delta_{c(j \to i) \omega_j, s_i c(j \to i)
  \omega_j}\,, 
\\[.1in]
\Delta_{j \to \overline i} 
= \Delta_{s_i c(j \to i) \omega_j, c(j \to i) \omega_j} \,.
\end{array}
\end{eqnarray*}
For a given $i$, we say that each minor of the form $\Delta_{j \to i}$
(resp. $\Delta_{j \to \overline i}$) is an $i$-\emph{indicator}
(resp. $\overline i$-indicator).

\begin{theorem} 
\label{th:GK general}
Let $C$ be a collection of $2r$ generalized minors that contains,
for every $i \in [1,r]$, an $i$-indicator and an $\overline i$-indicator.
Then, for an element $x \in G_{\geq 0}\,$, the following are
equivalent:
\begin{itemize}
\item[(g)]
$x$ is oscillatory; 
\item[(h)]
$\Delta (x) > 0$ for any $\Delta \in C$;
\item[(i)]
$x$ does not belong to a proper parabolic subgroup of $G$ containing 
$B$ or~$B_-$.  
\end{itemize}
\end{theorem}

Note that the equivalence ${\rm(g)}\Longleftrightarrow{\rm(h)}$ in
Theorem~\ref{th:GK general} generalizes Theorem~\ref{th:GK classical}.
Indeed, for $G = SL_n$ and the standard 
numbering of fundamental weights, one checks that 
$x_{i,i+1} = \Delta_{1 \to i}$ and $x_{i+1,i} = \Delta_{1 \to
  \overline i}\,$.  
Thus the set $C$ consisting of these matrix entries
satisfies the condition of Theorem~\ref{th:GK general}. 

Our last main result is a generalization of 
Theorem~\ref{th:GK bound} to all classical groups.

\begin{theorem} 
\label{th:GK2 general} 
For any given $G$, there exists a positive integer $m$ 
with the following property: an element $x \in G_{\geq 0}$ 
is oscillatory if and only if 
$x^m \in G_{>0}\,$.
A positive integer $m$ has this property if and only if 
for any permutation $\ii = (i_1, \dots, i_r)$ of indices $1, \dots, r$, 
the concatenation of $m$ copies of $\ii$ has a reduced word for
$\wnot$ as a subword.  
\end{theorem} 

Let $m(G)$ denote the smallest positive integer~$m$ that has the
property described in Theorem~\ref{th:GK2 general}.

\begin{theorem} 
\label{th:GK2 general concrete}
For a simple group~$G$,  
the value of $m(G)$ is given by the table 
\begin{center}
\begin{tabular}{|c|c|c|c|c|c|c|c|c|c|}
\hline
Type & $A_r$ & $B_r$ or $C_r$ & $D_r\,$, $r$ even
& $D_r\,$, $r$ odd  & $E_6$ & $E_7$ & $E_8$ & $F_4$ & $G_2$ \\
\hline
$m(G)$ & $r$ & $r$ & $r-1$
& $r$ & $8$ & $9$ & $15$ & $6$ & $3$ \\
\hline
\end{tabular} \,.
\end{center}
\end{theorem} 

The remaining sections contain the proofs of Theorems~\ref{th:loewner
  general}--\ref{th:GK2 general concrete}.  

\section{Proof of Theorem~\ref{th:TNN-TP general}}
\label{sec:TNN-TP general} 


The group $G$ has two \emph{Bruhat decompositions}, 
with respect to opposite Borel subgroups $B$ and $B_-\,$:
$$G = \bigcup_{u \in W} B u B = \bigcup_{v \in W} B_- v B_-  \ . $$
The \emph{double Bruhat cells}~$G^{u,v}$ are defined by 
$G^{u,v} = B u B  \cap B_- v B_- \,$.

Let $H_{>0}$ be the subgroup of~$H$ generated by the elements
$t^{h_i}$ for any $t > 0$ and $i \in [1,r]$; equivalently,
$H_{>0}$ consists of all $a\in H$ such that
$a^\gamma >  0$ for any weight $\gamma \in P$.  
Following G.~Lusztig, let us define the set~$G_{\geq 0}$ 
as the multiplicative monoid in $G$ generated by 
$H_{>0}$ and the elements $x_i (t)$ and $x_{\overline i} (t)$,  
for $i\in [1,r]$ and $t > 0$. 
In other words, we use condition (c) of Theorem~\ref{th:loewner
  general} as the interim definition of~$G_{\geq 0}$.  

The set $G_{\geq 0}$ is the disjoint union of \emph{totally positive
  varieties} $G^{u,v}_{> 0}$ 
defined by
$$G^{u,v}_{> 0} = G_{\geq 0} \cap G^{u,v} \, .$$

We denote $[\overline 1, \overline r] = \{\overline 1, \ldots,
\overline r\}$. 
For any sequence $\ii= (i_1, \ldots, i_m)$ of indices 
from the alphabet $[1,r] \cup [\overline 1, \overline r]$,
let us define the map $x_\ii: H \times \CC^m \to G$ by
\begin{equation}
x_\ii (a; t_1, \ldots, t_m) = a\, x_{i_1} (t_1) \cdots x_{i_m} (t_m) \, .
\end{equation} 
By definition, an element $x\in G_{\geq 0}$  
can be represented as $x=x_\ii (a; t_1, \ldots, t_m)$,
for some sequence $\ii$, with all the $t_k$
positive and $a\in H_{>0}\,$. 


A \emph{double reduced word} for the elements $u,v\in W$
is a reduced word for an element $(u,v)$ 
of the Coxeter group $W \times W$.
To avoid confusion, we will use the indices 
$\overline 1, \overline 2, \ldots, \overline r$ for the simple reflections 
in the first copy of $W$, and $1, 2, \ldots, r$ for the second copy.
A double reduced word for $(u,v)$ is nothing but a shuffle of a reduced
word for $u$ written in the alphabet 
$[\overline 1, \overline r]$ and a reduced word for $v$ written in 
the alphabet $[1,r]$.
We denote the set of double reduced words for $(u,v)$ by $R(u,v)$.

The \emph{weak order} is the partial order  on~$W$ defined as follows:
$u'\preceq u$ stands for $\l(u)=\l(u')+\l({u'}^{-1}u)$. 
(In other words, a reduced word for $u'$ can be extended on the right
to form a reduced word for~$u$.) 
We note that $w\preceq \wnot$ for any~$w \in W$.

The following lemma provides alternative
descriptions of the totally positive varieties~$G^{u,v}_{> 0}\,$. 

\begin{lemma}
\label{lem:c'''}
For an element $x\in G^{u,v}$, the following 
conditions are equivalent:
\begin{itemize}
\item[(${\rm c}'$)]
$x \in G^{u,v}_{> 0}\,$; 
\item[(${\rm c}''$)]
for some (equivalently, any) double reduced word 
$\ii \in R(u,v)$, we have $x = x_\ii (a;t_1, \dots, t_m)$ 
with $a \in H_{>0}$ and $t_1, \dots, t_m > 0$;
\item[(${\rm c}'''$)]
$\Delta_{u' \omega_i, v' \omega_i} (x) > 0$
for all $i\in [1,r]$ and all $u' \preceq u$, $v' \preceq v^{-1}$.
\end{itemize}
\end{lemma}

\proof
See \cite[Theorems~1.3 and 1.11]{FZ}. 
(The equivalence $({\rm c}')\Longleftrightarrow ({\rm c}'')$ was
essentially established in~\cite{lusztig}.)  
\endproof

Now everything is ready for the proof of Theorem~\ref{th:TNN-TP general}.
The implication $({\rm f})\Longrightarrow ({\rm d})$ is a special case of 
$({\rm c}')\Longrightarrow ({\rm c}''')$, while $({\rm d})\Longrightarrow ({\rm e})$ is 
trivial. 
Finally, to show that $({\rm e})\Longleftrightarrow ({\rm f})$, 
it suffices to note that
\[
G^{\wnot, \wnot} 
= \wnot G_0 \,\cap \, G_0 \wnot 
= \{x \in G: \Delta_{\wnot\omega_i,\omega_i}(x)\neq 0, \, 
\Delta_{\omega_i,\wnot\omega_i}(x)\neq 0 \, \, (i \in [1,r])\}
\] 
(cf.~\cite[Corollary~2.5]{FZ} or \cite[Proposition~4.1]{FZcells}). 
\qed


\section{Proof of Theorem~\ref{th:loewner general}}


\subsection{Proof of $({\rm b}) \Rightarrow ({\rm a})$}

This is obvious since all generalized minors are continuous functions
on~$G$. 

\subsection{Proof of $({\rm c}) \Rightarrow ({\rm b})$}

In view of Lemma~\ref{lem:c'''},
it suffices to show that the closure of $G^{\wnot, \wnot}_{> 0}$ 
contains all totally positive varieties $G^{u, v}_{> 0}\,$. 
Suppose $x \in G^{u,v}_{> 0}$ for some $u$ and $v$. 
Take any $\ii \in R(u,v)$ and write $x = x_\ii (a;t_1, \dots, t_m)$ 
as in 
(${\rm c}''$). 
Choose a word $\jj \in R(\wnot, \wnot)$
that has $\ii$ as an initial segment. 
Then
\begin{equation}
x = x_\ii (a;t_1, \dots, t_m) = 
\lim_{t \to +0} x_\jj (a;t_1, \dots, t_m, t, \dots, t) \,, 
\end{equation}
and $({\rm b})$ follows. 

\subsection{Proof of $({\rm a}) \Rightarrow ({\rm c})$}

Suppose that $x \in G^{u,v}$ satisfies $({\rm a})$.
It suffices to check condition $({\rm c}''')$ in Lemma~\ref{lem:c'''}.
Let $\Sigma (x)$ denote the set of all pairs
$(u', v') \in W \times W$ such that 
$\Delta_{u' \omega_i, v' \omega_i} (x) > 0$ for all $i$.
Our aim is to show that
\begin{equation}
\label{eq:aim}
(u',v') \in \Sigma (x) \ \ \text{for} \ \ u' \preceq u, \, 
v' \preceq v^{-1} \ .
\end{equation}
As a first step we notice that 
\begin{equation}
\label{eq:Schubert cell inequalities}
(u, e), (e, v^{-1}) \in \Sigma (x) \ ;
\end{equation}
this follows from the well-known fact that 
$\Delta_{u \omega_i, \omega_i}$ vanishes nowhere on 
the Bruhat cell $B u B$; see, e.g., \cite[Lemma~3.4]{FZcells}.

We shall write $u' \to u''$ if $u'' = u' s_i$ for some $i$, and 
$\l (u'') = \l (u') + 1$.  
In view of (\ref{eq:Schubert cell inequalities}), 
the desired inclusions (\ref{eq:aim}) are consequences of the
following statements:
\begin{eqnarray}
\label{eq:stat1}
\begin{array}{l}
\text{if $u' \to u''$ and $(u'', e)\in\Sigma (x)$,
then $(u', e)\in\Sigma (x)$;}\\[.1in]
\text{if $u' \to u''$ and $(e, u'')\in\Sigma (x)$,
then $(e, u')\in\Sigma (x)$;}
\end{array}
\end{eqnarray}
\begin{equation}
\label{eq:stat2}
\text{if $u' \to u''$, $v' \to v''$, $(u', v'')\in\Sigma (x)$, 
$(u'', v')\in\Sigma (x)$, 
then $(u'', v'') \in \Sigma(x)$.}
\end{equation} 
Our proof of both (\ref{eq:stat1}) and (\ref{eq:stat2}) 
relies on the following identity 
\cite[Theorem~1.17]{FZ}:
\begin{equation}
\label{eq:minors-Dodgson}
\Delta_{u'\omega_i, v' \omega_i} \Delta_{u'' \omega_i, v'' \omega_i}
= \Delta_{u' \omega_i, v'' \omega_i} \Delta_{u'' \omega_i, v' \omega_i}
+ \prod_{j \neq i} \Delta_{u' \omega_j, v' \omega_j}^{- a_{ji}} \,,
\end{equation}
whenever $u' \to u'' = u' s_i$ and $v' \to v'' = v' s_i$;
here the numbers $a_{ji}$ are the entries of the Cartan matrix of $G$. 

To prove (\ref{eq:stat1}), suppose that $u' \to u'' = u' s_i$ and 
$(u'', e) \in \Sigma (x)$. 
Now specialize (\ref{eq:minors-Dodgson}) at $v' = e$  and
evaluate both sides at~$x$.  
Using the fact that $u' \omega_j = u'' \omega_j$ for $j \neq i$, we see that
the second summand on the right-hand side is strictly positive.
Since all generalized minors of $x$ are nonnegative, we conclude 
that both factors on the left-hand side are positive.
In particular, $\Delta_{u'\omega_i, \omega_i} (x) > 0$, 
i.e., $(u', e) \in \Sigma (x)$, as desired. 
The second part of (\ref{eq:stat1}) is proved in the same way. 

To prove (\ref{eq:stat2}), suppose that $u' \to u'' = u' s_i$
and $v' \to v'' = v' s_j$, and both $(u', v'')$ 
and $(u'', v')$ belong to $\Sigma (x)$. 
We need to show that $\Delta_{u'' \omega_k, v'' \omega_k} (x) > 0$
for all $k$.
If $k \neq i$, then $u''\omega_k = u' \omega_k$ and we are done since
$(u', v'') \in \Sigma (x)$.
The case $k \neq j$ is treated in the same way. 
It thus remains to consider the case $k = j = i$. 
But then in (\ref{eq:minors-Dodgson}),
the first summand on the right (evaluated at $x$) is positive, implying
$\Delta_{u'' \omega_i, v'' \omega_i} (x) > 0$, as desired.
This completes the proof of Theorem~\ref{th:loewner general}.
\endproof

\section{Proof of Theorem~\ref{th:GK general}}

\subsection{Proof of $({\rm g}) \Rightarrow ({\rm i})$}

Since total positivity is described by condition $({\rm f})$ in
Theorem~\ref{th:TNN-TP general}, 
it suffices to show that every proper parabolic subgroup of $G$ containing 
$B$ or~$B_-$ has empty intersection with the open double Bruhat cell
$G^{\wnot, \wnot}$.  
The latter follows at once from the well known description of maximal
proper parabolic subgroups containing  
$B$ or~$B_-\,$: they are the subgroups $P_1, \dots, P_r$ and
$P_{\overline 1}, \dots, P_{\overline r}$ 
given by 
\begin{equation}
\label{eq:max parabolics}
P_i = \bigcup_{i \notin {\rm Supp} (v)} B_- v B_-\ ,
\quad P_{\overline i} = 
\bigcup_{i \notin {\rm Supp} (u)} B u B \ ,
\end{equation} 
where ${\rm Supp} (w)$ denotes the set of indices that occur in some (equivalently, any) reduced word for~$w \in W$. 

\subsection{Proof of $({\rm i}) \Rightarrow ({\rm g})$}

Consider the monoid $\mathcal{H}$ 
whose generators $T_1,\dots,T_r$ are subject to relations 
\[
\begin{array}{rcl}
T_i^2 &=& T_i \,;\\[.1in]
\underbrace{T_iT_jT_i\cdots}_{m_{ij}} &
=& \underbrace{T_jT_iT_j\cdots}_{m_{ij}} \ \ (i \neq j) \, ;
\end{array}
\]
here $m_{ij}$ is the order of $s_i s_j$ in $W$. 
A~well known theorem of Tits on reduced words (see, e.g.,
\cite[II,\S3C]{brown}) has the following implications. 
First, if $(i_1, \dots, i_m) \in R(w)$, then the product $T_{i_1}
\cdots T_{i_m}$ 
only depends on $w$ and so can be unambiguously denoted by~$T_w\,$. 
Second, the correspondence $w \mapsto T_w$ is a bijection between
$W$ and~$\mathcal{H}$. 
Finally, we have the following criterion for determining when a
product of generators is equal to $T_{\wnot}$. 

\begin{lemma}
\label{lem:hecke2}
For a word $(i_1,\dots,i_N)$ in the alphabet $[1,r]$,
we have $T_{i_1} \cdots T_{i_N}\!=\!T_{\wnot}$ if and only if this word
has a reduced word for~$\wnot$ as a subword. 
\end{lemma} 


The relevance of $\mathcal{H}$ to our problem is clear from the
following lemma.

\begin{lemma}
\label{lem:hecke}
For any $x \in G^{u,v}_{>0}$ and $y \in G^{u',v'}_{>0}$,
we have $xy \in G^{u'',v''}_{>0}$, 
where the elements $u''$ and $v''$ are given by $T_{u''}=T_u T_{u'}$
and $T_{v''}=T_v T_{v'}$. 
\end{lemma} 

\proof
Follows from condition $({\rm c}'')$ of
Lemma~\ref{lem:c'''}, together with the commutation relations among the
elementary factors $x_i(t)$ and $x_{\bar i}(t)$, as given in
\cite[Theorem~3.1]{BZ} and \cite[Section~2.2]{FZ}. 
\endproof

By Lemma~\ref{lem:hecke} and condition $({\rm f})$ of
Theorem~\ref{th:TNN-TP general},  
for any $x \in G^{u,v}_{>0}$ and any positive integer~$m$, we have
\begin{equation}
\label{eq:osc exponent}
x^m \in G_{>0} \Leftrightarrow T_u^m = T_v^m = T_{\wnot} \ .
\end{equation}
Suppose that a totally nonnegative element $x$ satisfies
condition~$({\rm i})$. 
By (\ref{eq:max parabolics}), $x \in G^{u,v}_{>0}$ for some elements
$u, v \in W$ such that  
${\rm Supp} (u) = {\rm Supp} (v) = [1,r]$. 
We need to show that $x$ is oscillatory.
In view of (\ref{eq:osc exponent}), this means 
that $T_u^m = T_v^m = T_{\wnot}$ for sufficiently large~$m$. 
The latter is clear from Lemma~\ref{lem:hecke2}: just
take~$m=\l(\wnot)$.

\subsection{Proof of $({\rm h}) \Leftrightarrow ({\rm i})$}

This equivalence can be restated as follows. 

\begin{lemma}
\label{lem:parabolic indicators}
Let $i \in [1,r]$, and let $\Delta$ be an $i$-indicator (resp.\ 
$\overline i$-indicator). 
Then $\Delta$ vanishes on $P_i$ (resp.\ $P_{\overline i}$),
and $\Delta (x) > 0$ for any $x \in G_{\geq 0}$ outside
$P_i$ (resp.~$P_{\overline i}$). 
\end{lemma}

\proof
It is enough to consider $i$-indicators, the case of 
$\overline i$-indicators being totally similar. 
Changing if necessary the numeration of fundamental weights,
we can assume without loss of generality that $i = 1$, and 
$$\Delta = \Delta_{j \to 1} = 
\Delta_{u \omega_j, s_1 u \omega_j} \, ,$$
where $u = c(j\to 1) = s_2 \cdots s_j\,$, with nonzero Cartan matrix 
entries $a_{k,k+1}$ for $k = 1, \dots, j-1$. 

First let us show that $\Delta (x) = 0$ for $x \in P_1$. 
We will denote by $x^T$ the ``transpose'' of~$x$;
more precisely, $x\mapsto x^T$ is the anti-automorphism of~$G$ defined
by
\begin{equation*}
\label{eq:T}
a^T = a \quad (a \in H) \ , \quad x_i (t)^T = x_{\overline i} (t) \ , 
\quad x_{\overline i} (t)^T = x_i (t) \,. 
\end{equation*} 
As in~\cite{FZ}, we will use the notation 
$\Delta^{\omega_i}=\Delta_{\omega_i, \omega_i}$
for the $i$th ``principal minor.'' 
Using \cite[(1.10), (2.25)]{FZ}, we obtain:
\[
\begin{array}{r}
\Delta (x) 
= \Delta_{u\omega_j,s_1u\omega_j}(x)
= \Delta_{s_1u\omega_j,u\omega_j}(x^T)
= \Delta^{\omega_j}({\overline{s_1u}}^{-1}x^T\overline{u})\\[.1in]
= \Delta^{\omega_j}({\overline{u^{-1}s_1u}}^{-1}\overline{u^{-1}}x^T\overline{u})
= \Delta_{u^{-1} s_1 u \omega_j, \omega_j} 
  (\overline {u^{-1}} x^T \overline u) \, .
\end{array}
\]
Observe that 
$\overline {u^{-1}} x^T \overline u \in P_{\overline 1}$ 
for any $x \in P_1$ (since all three factors belong to~$P_{\overline
  1}$). 
It remains to prove that $\Delta_{u^{-1} s_1 u \omega_j, \omega_j}$
vanishes on $P_{\overline 1}$. 
To see this we use the following description of $P_{\overline 1}$ 
equivalent to (\ref{eq:max parabolics}): 
$P_{\overline 1} = \pi^{-1} (X_{\wnot'})$, where 
$\pi$ is the projection of $G$ onto the flag variety $G/B$, 
the element $\wnot' \in W$ is the longest element of the parabolic
subgroup generated by $s_2, \dots, s_r$, and $X_w$ is the Schubert
variety corresponding to $w$ (i.e., the closure of the Schubert cell
$(B wB)/B$). 
Our claim that $\Delta_{u^{-1} s_1 u \omega_j, \omega_j}$
vanishes on $P_{\overline 1}$ now follows from the fact that
$1 \in {\rm Supp}(u^{-1} s_1 u)$, which means that $u^{-1} s_1 u$
is \emph{not} smaller than or equal to $\wnot'$ in the Bruhat order
(cf., e.g., \cite[Lemma 3.4]{FZcells};
in the notation of~\cite{FZcells},
$\Delta_{\gamma,\omega_j}(x)=p_\gamma(\pi(x))$). 

To complete the proof of Lemma~\ref{lem:parabolic indicators}
and Theorem~\ref{th:GK general}, it remains to show that
$\Delta_{j \to 1} (x) > 0$ for any element $x \in G_{\geq 0}$ 
not belonging to $P_1$.  
We proceed by induction on $j$. 
Let us first consider the case $j = 1$ when we need to show that
$\Delta_{\omega_1, s_1 \omega_1} (x) > 0$. 
Since 
$\Delta_{\omega_1, s_1 \omega_1}(b_- x) = 
\Delta^{\omega_1} (b_-) \Delta_{\omega_1, s_1 \omega_1}(x)$
for any $b_- \in B_-$, we can assume without loss of generality 
that $x$ has the form
$$x = x_{i_1} (t_1) \cdots x_{i_m}(t_m)$$
for some sequence of (unbarred) indices $i_1, \dots, i_m$ 
and some positive numbers $t_1, \dots, t_m$. 
The condition $x \notin P_1$ means that at least one of the indices
$i_k$ is equal to $1$; let $k$ be the maximal index such that $i_k = 1$.
Using the fact that 
$\overline {s_1}^{\ -1} x_i (t) \overline {s_1} \in N$ 
for any $i \neq 1$, and the commutation relation \cite[(2.13)]{FZ},
we conclude that
\[
\begin{array}{r}
\Delta_{\omega_1, s_1 \omega_1}(x) = 
\Delta^{\omega_1} (x_{i_1} (t_1) \cdots x_{i_k}(t_k) \overline {s_1} \cdot
(\overline {s_1}^{\ -1} x_{i_{k+1}} (t_{k+1}) \cdots x_{i_m}(t_m)
\overline {s_1})) \\[.1in]
=
\Delta^{\omega_1} (x_{i_1} (t_1) \cdots x_{i_k}(t_k) \overline {s_1}) 
= \Delta^{\omega_1} (x_{i_1} (t_1) \cdots x_{i_{k-1}}(t_{k-1})
x_{\overline 1}(t_k^{-1}) t_k^{h_1})\, .
\end{array}
\]
Since the element $x' = x_{i_1} (t_1) \cdots x_{i_{k-1}}(t_{k-1})
x_{\overline 1}(t_k^{-1}) t_k^{h_1}$ is totally nonnegative, 
and any principal minor is positive on~$G_{\geq 0}$
(see \cite[Corollary~2.5 and Proposition~2.29]{FZ}), 
we conclude that $\Delta^{\omega_1} (x') > 0$, as desired.

Now assume that $j \geq 2$, and that we already know 
that $\Delta_{j' \to 1}(x) > 0$ for $j' = 1, \dots, j-1$. 
Let us apply the identity (\ref{eq:minors-Dodgson}) for 
$i = j$, $u' = s_2 \cdots s_{j-1}$, and $v' = s_1 \cdots s_{j-1}$.
In our present notation, it takes the following form:
\begin{equation}
\label{eq:Dodgson indicator}
\Delta^{\omega_j} \Delta_{j \to 1} = 
\Delta_{u' \omega_j, v'' \omega_j} \Delta_{u'' \omega_j, v' \omega_j}
+ \prod_{j' > j} (\Delta^{\omega_{j'}})^{- a_{j'j}} \cdot
\prod_{j'= 1}^{j-1} \Delta_{j' \to 1}^{- a_{j'j}} \,; 
\end{equation}
here we used that $\Delta_{u\omega_j,v\omega_j}=\Delta^{\omega_j}$
whenever $u$ and $v$ belong to the parabolic subgroup of~$W$ generated
by all simple reflections except~$s_j\,$.  
By the inductive assumption, the second summand in the right-hand side
of~(\ref{eq:Dodgson indicator}) is positive at~$x$, while the first summand 
is nonnegative. 
It follows that $\Delta_{j \to 1}(x) > 0$, completing the proof.
\endproof

\section{Proof of Theorems~\ref{th:GK2 general} and \ref{th:GK2
    general concrete}} 

\subsection{Proof of Theorem~\ref{th:GK2 general}}

This is an immediate consequence of (\ref{eq:osc exponent}) and
Lemma~\ref{lem:hecke2}. 

\subsection{Proof of Theorem~\ref{th:GK2 general concrete}}
\label{sec:coxeter-elements}

We will need some basic facts about Coxeter elements in Weyl groups
(the proofs can be found in~\cite[Section~V.6]{bourbaki}). 
Recall that a \emph{Coxeter element} $c \in W$ is a product of simple
reflections $s_1, \dots, s_r$ taken in any order. 
All such elements are conjugate to each other and thus have the same order;
this order is called the \emph{Coxeter number} of $W$ and denoted by~$h$. 
Here are the statements we need:
\begin{itemize}
\item[(C1)] 
If $W$ is irreducible, then $h/2 = \l(\wnot)/r$. 
\item[(C2)]
If $\wnot = -1$ (i.e., $\wnot (\lambda) = -
\lambda$ for any weight $\lambda$), 
then $h$ is even, and $c^{h/2} = \wnot$ for any Coxeter element $c \in
W$. 
\end{itemize} 

Now suppose that $G$ is simple, so the Weyl group $W$ is irreducible.
Combining Theorem~\ref{th:GK2 general} with (C1)--(C2),
we conclude that $m(G) = h/2 = \l(\wnot)/r$ whenever $\wnot = -1$. 
According to the tables in~\cite{bourbaki}, this gives the desired
answer for $m(G)$ for all the types 
except $A_r$ ($r \geq 2$), $D_r$ ($r$ odd), and~$E_6\,$. 
Let us consider these remaining cases separately.

Throughout, we denote by $\ii=(i_1,\dots,i_r)$ a permutation of
indices~$1,\dots,r$. 
It will be convenient to use the notation $\ii^k$ for the
concatenation of $k$ copies of~$\ii$. 

\emph{Type $A_r\,$}. 
As usual, we identify $W$ with the symmetric group $S_{r+1}$; under this identification,
$s_i$ becomes the transposition of adjacent indices $i$ and $i+1$, and $\wnot (i) = r + 2 - i$
for $i \in [1,r+1]$. 
If $\ii=(1,\dots,r)$, then~$\ii^{r-1}$ does not
contain a reduced word for~$\wnot$, 
since any such reduced word must have a subword
$r,r\!-\!1,\dots,2,1$ (because $\wnot$ switches 1 and~$r+1$). 
For an arbitrary permutation~$\ii$ of $1,\dots,r$,
let us now consider the sequence~$\ii^r$.
We will form a subsequence~$\jj$ of~$\ii^r$  as follows. 
First, $\jj$ will include all $r$ entries of~$\ii^r$ which are equal
to~1. 
Between any two consecutive~1's, there is a~2;
let $\jj$ include all these~2's (there will be $r-1$ of them).
We then include in~$\jj$ the~3's that interlace these~2's
($r-2$ more entries), etc. 
It is straightforward to check that the subsequence~$\jj$ thus
obtained will be a reduced word for $\wnot\,$.
Thus $m(G) = r$, as claimed. 

\emph{Type $D_r$ ($\,r$ odd).}
In this case $h/2 = \l(\wnot)/r = r-1$. 
Using the standard combinatorial interpretation of~$D_r\,$, one checks
that $(s_1\cdots s_r)^{r-1}\neq\wnot\,$, and so $m(G) \geq r$. 
To prove the reverse inequality, consider the standard embedding of $W$ 
into the Coxeter group $\widetilde W$ of type $D_{r+1}$. 
We know that the Coxeter number $\tilde h$ of $\widetilde W$ is equal to~$2r$,
and the longest element $\tilde \wnot \in \widetilde W$ is equal to $-1$. 
Let $\ii=(i_1,\dots,i_r)$ be a permutation of $1,\dots,r$,
and denote $\tilde \ii=(i_1,\dots,i_r,r+1)$. 
Then ${\tilde \ii}^r$ is a reduced
word for $\tilde \wnot$, and therefore it
contains a reduced word for $\wnot \in W$ as a subword.
We conclude that $m(G) = r$, as desired. 

\begin{figure}[ht]
\setlength{\unitlength}{2pt} 
\begin{center}
\begin{picture}(80,25)(0,-5)
\thicklines

\put(0,0){\line(1,0){80}}
\put(40,0){\line(0,1){20}}

\put(0,0){\circle*{2.5}}
\put(20,0){\circle*{2.5}}
\put(40,0){\circle*{2.5}}
\put(60,0){\circle*{2.5}}
\put(80,0){\circle*{2.5}}
\put(40,20){\circle*{2.5}}

\put(0,-5){\makebox(0,0){$s_1$}}
\put(20,-5){\makebox(0,0){$s_2$}}
\put(40,-5){\makebox(0,0){$s_4$}}
\put(60,-5){\makebox(0,0){$s_5$}}
\put(80,-5){\makebox(0,0){$s_6$}}
\put(35,20){\makebox(0,0){$s_3$}}

\end{picture}

\end{center}
\caption{Generators of the Weyl group of type~$E_6$}
\label{fig:E6}

\end{figure}
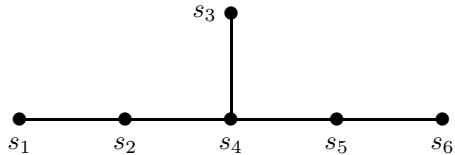

\emph{Type $E_6$}. 
The upper bound $m(G)\leq 8$ 
can be proved using the fact that $(s_1s_4s_6s_2s_3s_5)^6=1$
(in the notation of Figure~\ref{fig:E6}),
together with the following observation based on Lemma~\ref{lem:hecke2}:
if $(T_c)^k=T_{\wnot}$ for a Coxeter element~$c\in W$,
then $(T_{c'})^{k+1}=T_{\wnot}$ for any Coxeter element~$c'$ obtained
by taking a cyclic permutation of any reduced word for~$c$. 
The lower bound is proved by exhibiting a Coxeter element
(namely, $c=s_1s_2s_3s_4s_5s_6$) such that $(T_c)^7\neq
T_{\wnot}\,$. 
(The latter verification is due to H.~Derksen.) 
\endproof

\section*{Acknowledgements}

Harm Derksen contributed to
Section~\ref{sec:coxeter-elements} by first bringing 
the statements (C1)--(C2) to our attention,
and then by verifying the type~$E_6$ case of
Theorem~\ref{th:GK2 general concrete}.  
We are grateful to Harm for his input.

\end{document}